\documentstyle[12pt]{article} 
\def\addto#1#2{
\ifx\zone\undefined\let\zone=#1\def#1{\zone#2}\else
\ifx\ztwo\undefined\let\ztwo=#1\def#1{\ztwo#2}\else
\ifx\zthree\undefined\let\zthree=#1\def#1{\zthree#2}\else
\ifx\zfour\undefined\let\zfour=#1\def#1{\zfour#2}\else
\ifx\zfive\undefined\let\zfive=#1\def#1{\zfive#2}\else
\ifx\zsix\undefined\let\zsix=#1\def#1{\zsix#2}\else
\ifx\zsevn\undefined\let\zsevn=#1\def#1{\zsevn#2}\else
\ifx\zegt\undefined\let\zegt=#1\def#1{\zegt#2}\else
\ifx\znin\undefined\let\znin=#1\def#1{\znin#2}\else
\ifx\zten\undefined\let\zten=#1\def#1{\zten#2}\else
\fi\fi\fi\fi\fi\fi\fi\fi\fi\fi
}

\expandafter\chardef\csname pre amsym.def at\endcsname=\the\catcode`\@
\catcode`\@=11
\def\undefine#1{\let#1\undefined}
\def\newsymbol#1#2#3#4#5{\let\next@\relax
 \ifnum#2=\@ne\let\next@\msafam@\else
 \ifnum#2=\tw@\let\next@\msbfam@\fi\fi
 \mathchardef#1="#3\next@#4#5}
\def\mathhexbox@#1#2#3{\relax
 \ifmmode\mathpalette{}{\m@th\mathchar"#1#2#3}%
 \else\leavevmode\hbox{$\m@th\mathchar"#1#2#3$}\fi}
\def\hexnumber@#1{\ifcase#1 0\or 1\or 2\or 3\or 4\or 5\or 6\or 7\or 8\or
 9\or A\or B\or C\or D\or E\or F\fi}
\newfam\msbfam
\edef\msbfam@{\hexnumber@\msbfam}
\def\Bbb#1{\fam\msbfam\relax#1}
\@ifundefined{setboxz@h}{\def\setboxz@h{\setbox\z@\hbox}}{}  
\@ifundefined{wdz@}{\def\wdz@{\wd\z@}}{}                     
\def\widehat#1{\setboxz@h{$\m@th#1$}%
 \ifdim\wdz@>\tw@ em\mathaccent"0\msbfam@5B{#1}%
 \else\mathaccent"0362{#1}\fi}
\catcode`\@=\csname pre amsym.def at\endcsname
\newfont{\seveneufm}{eufm8 scaled 1000} 
\newfont{\egteufm}{eufm10 scaled 1000}  
\newfont{\nineufm}{eufm10 scaled 1095}    
\newfont{\teneufm}{eufm10 scaled 1200}  
\newfam\eufmfam
\def\xiieufm{\textfont\eufmfam=\teneufm
\scriptfont\eufmfam=\seveneufm
\def\got{\fam\eufmfam\teneufm}}
\def\xieufm{\textfont\eufmfam=\nineufm
\def\got{\fam\eufmfam\nineufm}}
\def\xeufm{\textfont\eufmfam=\egteufm
\def\got{\fam\eufmfam\egteufm}}
\newfont{\seveneusm}{eusm8 scaled 1000} 
\newfont{\egteusm}{eusm10 scaled 1000}  
\newfont{\nineusm}{eusm10 scaled 1095}    
\newfont{\teneusm}{eusm10 scaled 1200}  
\newfam\eusmfam
\def\xiieusm{\textfont\eusmfam=\teneusm
\scriptfont\eusmfam=\seveneusm
\def\skr{\fam\eusmfam\teneusm}}
\def\xieusm{\textfont\eusmfam=\nineusm
\def\skr{\fam\eusmfam\nineusm}}
\def\xeusm{\textfont\eusmfam=\egteusm
\def\skr{\fam\eusmfam\egteusm}}
\newfont{\sevenmsb}{msbm8 scaled 1000} 
\newfont{\egtmsb}{msbm10 scaled 1000}  
\newfont{\ninmsb}{msbm10 scaled 1095}    
\newfont{\tenmsb}{msbm10 scaled 1200}  
\def\xiimsb{\textfont\msbfam=\tenmsb
\scriptfont\msbfam=\sevenmsb}
\def\ximsb{\textfont\msbfam=\ninmsb}
\def\xmsb{\textfont\msbfam=\egtmsb}
\addto\normalsize{\xiieusm\xiieufm\xiimsb}  
\addto\small{\xieusm\xieufm\ximsb}  
\addto\footnotesize{\xeusm\xeufm\xmsb}  
\newtheorem{theorem}{Theorem}%
\newtheorem{lemma}[theorem]{Lemma}%
\newtheorem{definition}[theorem]{Definition}%
\newtheorem{corollary}[theorem]{Corollary}%
\newtheorem{proposition}[theorem]{Proposition}%
\newtheorem{property}[theorem]{Property}
\newcommand{\bte}{\begin{theorem}\ }%
\newcommand{\ete}{\end{theorem}}%
\newcommand{\ble}{\begin{lemma}\ }%
\newcommand{\ele}{\end{lemma}}%
\newcommand{\bdf}{\begin{definition}\ \rm}%
\newcommand{\edf}{\end{definition}}%
\newcommand{\bcor}{\begin{corollary}\ }%
\newcommand{\ecor}{\end{corollary}}%
\newcommand{\bpro}{\begin{proposition}\ }%
\newcommand{\epro}{\end{proposition}}%
\newcommand{\bprp}{\begin{property}\ }%
\newcommand{\eprp}{\end{property}}%
\newcommand{\ben}{\begin{enumerate}}%
\newcommand{\een}{\end{enumerate}}%
\newcommand{\bit}{\begin{itemize}}%
\newcommand{\eit}{\end{itemize}}%
\newcommand{\bay}{\begin{array}}%
\newcommand{\eay}{\end{array}}%
\newcommand{\ang} [1]{\langle #1\rangle}%
\newcommand{\ans} [1]{\{\hspace{0.2mm}#1\hspace{0.2mm}\}}%
\newcommand {\mem}{$\hspace{0.5pt}\msur\in\msur$} 
\newcommand{\res}{{\mathbin{\kern1pt|\kern 0.5pt}}}%
\newcommand{\we}{{\mathbin{\kern 1.3pt ^\wedge}}}%
\newcommand{\ain}{{\mathbin{\hspace{1pt}{\mathord{\in}}^\ast\hspace{1pt}}}}
\newcommand{\los}{\L o\'s}
\newcommand{\dom}{{\rm dom}\,}%
\newcommand{\ult}{{\rm Ult}}%
\newcommand{\Ord}{{\rm Ord}}%
\newcommand{\dd}[2]{\mbox{$#1\msur$-#2}}%
\newcommand{\its}{\vspace{-1mm}}%
\newcommand{\qed}{\hfill{$\msur\Box\msur$}}%
\newcommand{\lra}{\longrightarrow} %
\newcommand{\ti}{\times}%
\newcommand{\noi}{\noindent}%
\newcommand{\ccc}{\hspace{0.5pt}}%
\newcommand{\scri}[1]{\ccc\skr #1\ccc}%
\newcommand{\cN}{\mathord{\scri N}}%
\newcommand{\cU}{\mathord{\scri U}}%
\newcommand{\cF}{\mathord{\scri F}}%
\newcommand{\cL}{\mathord{\scri L}}%

\newcommand{\Ba}{{\bf a}}
\newcommand{\bbb}{\hspace{0.5pt}}
\newcommand{\bbL}{\mathord{\bbb\Bbb L\bbb}} %
\newcommand{\bbV}{\mathord{\bbb\Bbb V\bbb}} 
 
\newcommand{\bbP}{\mathord{\bbb\Bbb P\bbb}} 
\newcommand{\gM}{{\hspace{0.5pt}{\got M}\hspace{0.5pt}}} 
\newcommand{\gN}{{\hspace{0.5pt}{\got N}\hspace{0.5pt}}}
\newcommand{\ide}{\mathord{\hspace{0.5pt}{\got i}\hspace{1pt}}}
\newcommand{\ra}{\mathord{\hspace{0.5pt}{\got r}\hspace{1pt}}}
\newcommand{\al}{\alpha} %
\newcommand{\ba}{\beta}%
\newcommand{\ga}{\gamma}%
\newcommand{\om}{\omega} %
\newcommand{\vit}{\mathbin{{\hspace{0.5pt}\tt vit\hspace{0.5pt}}}}
\newcommand{\sq}{\subseteq}%
\newcommand{\dm}{$$}%

\newcommand{\cir}{\circ}
\newcommand{\proof}{\noi{\bf Proof} \ }%
\newcommand{\wh}{\widehat}
\newcommand{\un}{\underline}
\newcommand{\msur}{\hspace{-1mm}}%
\begin{document}  
\title{On a Spector ultrapower of the Solovay model\thanks
{Research supported by the Netherlands Organization for 
Scientific Research NWO under grant PGS 22 262}
}
\author{
Vladimir Kanovei\thanks{Moscow Transport Engineering Institute,  
\ {\tt kanovei@sci.math.msu.su}} 
\and 
Michiel van Lambalgen\thanks{University of Amsterdam, \ 
{\tt michiell@fwi.uva.nl}}
}
\date{\today}
\maketitle 

\thispagestyle{empty}

\begin{abstract}
We prove that a Spector--like ultrapower extension $\gN$ of a 
countable Solovay model $\gM$ (where all sets of reals are 
Lebesgue measurable) is equal to the set of all sets 
constructible from reals in a generic extension $\gM[\al]$ 
where $\al$ is a random real over $\gM.$ The proof involves 
an almost everywhere uniformization theorem in the Solovay 
model.
\end{abstract}


\newpage

\normalsize

\section*{Introduction}

Let $\cU$ be an ultrafilter in a transitive model $\gM$ of 
${\bf ZF}.$ Assume that an ultrapower of $\gM$ via $\cU$ is to be 
defined. The first problem we meet is that $\cU$ may not be an 
ultrafilter in the universe because not all subsets of the index 
set belong to $\gM$.

We can, of course, extend $\cU$ to a true ultrafilter, say 
$\cU',$ but this may cause additional trouble. Indeed, if $\cU$ 
is a special ultrafilter in $\gM$ certain properties of which 
were expected to be exploit, then most probably these 
properties do not transfer to $\cU';$ assume for instance that 
$\cU$ is countably complete in $\gM$ and $\gM$ itself is 
countable. Therefore, it is better to keep $\cU$ rather than 
any of its extensions in the universe, as the ultrafilter. 

If $\gM$ models ${\bf ZFC},$ the problem can be solved by taking 
the {\it inner\/} ultrapower. In other words, we consider only 
those functions $f:I\;\lra\;\gM$ (where $I\in\gM$ is the 
carrier of $\cU$) which belong to $\gM$ rather than {\it all\/} 
functions $f\in\gM^I,$ to define the ultrapower. This version, 
however, depends on the axiom of choice in $\gM;$ otherwise the 
proofs of the basic facts about ultrapowers (e.\ g. \los'   
theorem) will not work. 

The ``choiceless'' case can be handled by a sophisticated 
construction of {\sc Spector}~[1991], which is based on ideas 
from both forcing and the ultrapower technique. 
As presented in {\sc Kanovei} and {\sc van Lambalgen} 
[1994], this construction  
proceeds as follows.  One has to add to the 
family of functions $\cF_0=\gM^I\cap \gM$ a number of new 
functions $f\in\gM^I\,,$ $f\not\in\gM\,,$ which are intended 
to be choice functions whenever we need such in the ultrapower 
construction. 
 
In this paper, we consider a very interesting choiceless 
case: $\gM$ is a Solovay model of ${\bf ZF}$ plus the 
principle of dependent choice, in which all sets of reals 
are Lebesque measurable, and the ultrafilter $\cL$ on the set 
$I$ of Vitali degrees of reals in $\gM,$ generated by sets of 
positive measure. 

\newpage

\section{On a.e. uniformization in the Solovay model}
\label{sun}

In this section, we recall the uniformization properties  
in a Solovay model. Thus let $\gM$ be a countable transitive 
Solovay model for Dependent Choices plus ``all sets are Lebesgue 
measurable'', as it is defined in {\sc Solovay}~[1970], -- 
the {\it ground model\/}. The following known 
properties of such a model will be of particular interest below. 

\bprp
\label{sol1}
{\rm [True in $\gM$]} \\[1pt] 
$\bbV=\bbL({\rm reals})\,;$ in particular 
every set is real--ordinal--definable.\qed
\eprp

To state the second property, we need to introduce some notation. 

Let $\cN=\om^\om$ denote the Baire space, the elements of which 
will be referred to as {\it real numbers\/} or {\it reals\/}.. 

Let $P$ be a set of pairs such that $\dom P\sq\cN$ (for 
instance, $P\sq\cN^2$). We say that a function $f$ defined on $\cN$ 
{\it uniformizes\/} $P$ {\it a.e.\/} (almost everywhere) iff the set
\dm
\ans{\al\in\dom P:\ang{\al,f(\al)}\not\in P}
\dm
has null measure. For example if the projection $\dom P$ is a 
set of null measure in $\cN$ then any $f$ uniformizes a.e. $P,$ but 
this case is not interesting. The interesting case is the case when 
$\dom P$ is a set of full measure, and then $f$ a.e. uniformizes $P$ 
iff for almost all $\al,$ $\ang{\al,f(\al)}\in P_\al$. 

\bprp
\label{sol2}
{\rm [True in $\gM$]}\\[1pt] 
Any set\/ $P\in \gM\,,\;\,P\sq \cN^2,$ can be uniformized 
a.e. by a Borel function. {\rm (This implies the Lebesgue 
measurability of all sets of reals, which is known to be 
true in $\gM$ independently.)}\qed
\eprp

This property can be expanded (with the loss of the condition 
that $f$ is Borel) on the sets $P$ which do not necessarily 
satisfy $\dom P\sq \cN$.

\bte
\label{un}
In\/ $\gM,$ any set\/ $P$ with\/ $\dom P\sq \gM$ 
admits an a.e. uniformisation.  
\ete
\proof  Let $P$ be an arbitrary set of pairs such that 
$\dom P\sq \cN$ in $\gM.$ Property~\ref{sol1} implies the 
existence of a function $D:(\Ord\cap\gM)\ti 
(\cN\cap\gM)\;\;\hbox{onto}\;\;\gM$ which is \mem-definable 
in $\gM.$ 

{\it We argue in $\gM$.} 
Let, for $\al\in\cN,$ $\xi(\al)$ denote the least ordinal $\xi$ 
such that
\dm
\exists\,\ga\in\cN\;[\,\ang{\al,D(\xi,\ga)}\in P\,]\,.
\dm 
(It follows from the choice of $D$ that $\xi(\al)$ is well 
defined for all $\al\in\cN.$) It remains to apply 
Property~\ref{sol2} to the set 
$P'=\ans{\ang{\al,\ga}\in \cN^2:\ang{\al,D(\xi(\al),\ga)}\in P}$.
\qed


\section{The functions to get the Spector ultrapower}

We use a certain ultrafilter over the set of Vitali degrees of 
reals in $\gM,$ the initial Solovay model, to define the 
ultrapower. 

Let, for $\al,\,\al'\in\cN\,,$ $\al\vit\al'$ if and only if 
$\exists\,m\;\forall\,k\>m\;(\al(k)=\al'(k)),$  (the 
{\it Vitali equivalence\/}).\its
\bit
\item 
For $\al\in\cN,$ we set $\un\al=\ans{\al':\al'\vit\al},$ the   
{\it Vitali degree\/} of $\al$.\its

\item $\un\cN=\ans{\un\al:\al\in\cN}\,;$ $i,\,j$ denote elements 
of $\un\cN$.\its
\eit
As a rule, we shall use {\it underlined\/} characters $\un f,\;
\un F,\,...$ to denote functions defined on $\un\cN,$ while 
functions defined on $\cN$ itself will be denoted in the usual 
manner. 

Define, in $\gM,$ an ultrafilter $\cL$ 
over $\un\cN$ by: $\un X\sq\un\cN$ belongs to $\cL$ iff the 
set $X=\ans{\al\in\cN:\un\al\in\un X}$ has full Lebesgue measure. 
It is known (see e.g. {\sc van Lambalgen}~[1992], Theorem 2.3) 
that the measurability hypothesis implies that 
$\cL$ is \dd\kappa complete in $\gM$ for all cardinals $\kappa$ 
in $\gM$. 

One cannot hope to define a good \dd\cL ultrapower of $\gM$ using 
only functions from $\cF_0= \ans{\un f\in\gM:\dom \un f=\un\cN}$ 
as the 
base for the ultrapower. Indeed consider the identity function 
$\ide\in\gM$ defined by $\ide(i)=i$ for all $i\in\un\cN.$ Then 
$\ide(i)$ is nonempty for all $i\in\un\cN$ in $\gM,$ therefore to 
keep the usual properties of ultrapowers we need a 
function $\un f\in \cF_0$ such that $\un f(i)\in i$ for almost 
all $i\in\un\cN,$ but Vitali showed that such a choice function 
yields a nonmeasurable set. 

Thus at least we have to add to $\cF_0$ a new function $\un f,$ not 
an element of $\gM,$ which satisfies $\un f(i)\in i$ for almost 
all $i\in\un\cN.$ Actually it seems likely that we have to add 
a lot of new functions, to handle similar situations, including 
those functions the existence of which is somehow implied by the 
already added functions. A general way how to do this, extracted 
from the exposition in {\sc Spector}~[1991], was presented in 
{\sc Kanovei} and {\sc van Lambalgen} [1994]. However in the 
case of the Solovay model the 
a.e. uniformization theorem (Theorem~\ref{un}) allows to add 
essentially a single new function, corresponding to 
the \dd{\ide} case considered above. 

\subsubsection*{The generic choice function for the identity}

Here we introduce a function $\ra$ defined on $\cN\cap \gM$ and 
satisfying $\ra(i)\in i$ for all $i\in\cN\cap\gM.$ $\ra$ will be 
generic over $\gM$ for a suitable notion of forcing. 

The notion of forcing is introduced as follows. In $\gM,$ 
let $\bbP$ be the set of all functions $p$ defined on $\un\cN$ 
and satisfying $p(i)\sq i$ and $p(i)\not=\emptyset$ 
for all $i.$\footnote
{Or, equivalently, the collection of all sets $X\sq\cN$ 
which have a nonempty intersection with 
every Vitali degree. Perhaps this 
forcing is of separate interest. 
} 
(For example $\ide\in\bbP.$) We order $\bbP$ so 
that $p$ is stronger than $q$ iff $p(i)\sq q(i)$ for all $i.$ 
If $G\sq\bbP$ is \dd\bbP generic over $\gM$, $G$ defines a 
function $\ra$ by
\dm
\ra(i)=\hbox{the single element of}\;\;{\textstyle\bigcap
_{p\in G}}\;p(i)\,
\dm
for all $i\in\un\cN\cap\gM.$ Functions $\ra$ defined this way 
will be called \dd\bbP {\it generic over\/} $\gM.$ Let us fix 
such a function $\ra$ for the remainder of this paper.

\subsubsection*{The set of functions used to define the 
ultrapower}

We let $\cF$ be the set of all superpositions $f\cir \ra$ 
where\footnote
{To make things clear, $f\cir\ra(i)=f(\ra(i))$ for all $i$.}\   
$\ra$ is the generic function fixed above while $f\in\gM$ 
is an arbitrary function defined on $\cN\cap\gM.$ Notice that in 
particular any function $\un f\in\gM$ defined on $\un\cN\cap\gM$ 
is in $\cF:$ take $f(\al)=\un f(\un\al)$. 

To see that $\cF$ can be used successfully as the base of an 
ultrapower of $\gM,$ we have to check three fundamental conditions 
formulated in {\sc Kanovei} and {\sc van Lambalgen} [1994].

\bpro
\label{promes}
{\rm [Measurability]} \ Assume that\/ $E\in \gM$ and\/ 
$f_1,...,f_n\in \cF.$ Then the set\/ 
$\ans{i\in\un\cN\cap\gM:E(f_1(i),...,f_n(i))}$ 
belongs to $\gM$.
\epro
\proof By the definition 
of $\cF,$ it suffices to prove that $\ans{i:\ra(i)\in E}\in \gM$ 
for any set $E\in\gM,\;E\sq\cN.$ By the genericity of $\ra,$ it 
remains then to prove the following in $\gM:$ for any $p\in\bbP$ 
and any set $E\sq\cN,$ there exists a stronger condition $q$ 
such that, for any $i,$ either $q(i)\sq E$ or $q(i)\cap E=
\emptyset.$ But this is obvious.\qed 

\bcor
\label{mes}
Assume that\/ $V\in\gM,\;V\sq\cN$ is a set of null measure 
in\/ $\gM.$ Then, for \dd\cL almost all\/ $i,$ we have\/ 
$\ra(i)\not\in V$.
\ecor
\proof By the proposition, the set 
$I=\ans{i:\ra(i)\in V}$ belongs to $\gM.$ Suppose that, on the 
contrary, $I\in\cL.$ Then $A=\ans{\al:\un\al\in I}$ is a set of 
full measure. On the other hand, since $\ra(i)\in i,$ we have 
$A\sq\bigcup_{\ba\in V}\,\un\ba,$ where the right--hand side is 
a set of null measure because $V$ is such a set, 
contradiction.\qed

\bpro
\label{procho}
{\rm [Choice]} \ Let\/ $f_1,...,f_n\in \cF$ and\/ $W\in \gM.$ 
There exists  a function\/ $\un f\in \cF$ such that, for\/  
\dd\cL almost all $i\in\un\cN\cap\gM,$ it is true in\/ $\gM$ that
\dm
\exists\,x\;W(f_1(i),...,f_n(i),x)\,
\;\lra\;\,W(f_1(i),...,f_n(i),\un f(i))\,.
\dm
\epro
\proof This can be reduced to the following: given $W\in\gM,$ 
there exists a function $\un f\in\cF$ such that, for 
\dd\cL almost all $i\in\un\cN\cap\gM$, 
\dm
\exists\,x\;W(\ra(i),x)\,\;\lra\;\,W(\ra(i),\un f(i))
\eqno{(\ast)}
\dm
in $\gM.$ 

{\it We argue in $\gM$.} Choose $p\in \bbP.$ and 
let $p'(i)=\ans{\ba\in p(i):
\exists\,x\,W(\ba,x)},$ and $X=\ans{i:p'(i)\not=\emptyset}.$ 
If $X\not\in\cL$ then an arbitrary $\un f$ defined on $\un\cN$ 
will satisfy $(\ast),$ therefore it is assumed that $X\in\cL.$ 
Let 
\dm
q(i)=\left\{
\bay{ccl}
p'(i) & \hbox{iff} & i\in X\\[3mm]

p(i)  &            & \hbox{otherwise}
\eay
\right.
\dm
for all $i\in\un\cN;$ then $q\in\bbP$ is stronger than $p.$ 
Therefore, since $\ra$ is generic, one may assume that 
$\ra(i)\in q(i)$ for all $i.$ 

Furthermore, ${\rm DC}$ in the Solovay model $\gM$ implies that 
for every $i\in X$ the following is true: 
there exists a function $\phi$ defined on $q(i)$ and such that 
$W(\ba,\phi(\ba))$ for every $\ba\in q(i).$ Theorem~\ref{un} 
provides a function $\Phi$ such that for almost all $\al$ the 
following is true: the value $\Phi(\al,\ba)$ is defined and 
satisfies $W(
\ba,\Phi(\al,\ba))$ for all 
$\ba\in q(\un\al).$ Then, by Corollary~\ref{mes}, we have 
\dm
\hbox{for all}\;\;\ba\in q(\un{\ra(i)})\,,\;\;W(
\ba,\,\Phi(\ra(i),\ba)\,)
\dm
for almost all $i.$ However, $\un{\ra(i)}=i$ for all $i.$ 
Applying the assumption that $\ra(i)\in q(i)$ for all $i,$ we 
obtain $W(\ra(i),\,\Phi(\ra(i),\ra(i))\,)$ for almost all $i.$ 
Finally the function $\un f(i)=\Phi(\ra(i),\ra(i))$ is in 
$\cF$ by definition. 
\qed

\bpro
\label{proreg}
{\rm [Regularity]} \ For any\/ $\un f\in \cF$ there exists an 
ordinal\/ $\xi\in\gM$ such that for\/ \dd\cL almost all\/ $i,$ 
if\/  $\un f(i)$ is an ordinal then $\un f(i)=\xi$. 
\epro
\proof To prove this statement, assume that $\un f=f\cir \ra$ 
where $f\in\gM$ is a function defined on $\cN$ in $\gM$. 

{\it We argue in $\gM$.} Consider an arbitrary $p\in\bbP.$ We 
define a stronger condition $p'$ as follows. Let $i\in\un\cN.$ 
If there does not exist $\ba\in p(i)$ such that 
$f(\ba)$ is an ordinal, we put $p'(i)=p(i)$ and $\xi(i)=0.$ 
Otherwise, let $\xi(i)=\xi$ be the 
least ordinal $\xi$ such that $f(\ba)=\xi$ for some $\ba\in 
p(i).$ We set $p'(i)=\ans{\ba\in p(i):f(\ba)=\xi(i)}$.

Notice that $\xi(i)$ is an ordinal for all $i\in\un\cN.$ 
Therefore, since the ultrafilter $\cL$ is \dd\kappa  
complete in $\gM$ for all $\kappa,$ 
there exists a single ordinal $\xi\in\gM$ 
such that $\xi(i)=\xi$ for almost all $i$.

By genericity, we may assume that actually $\ra(i)\in
p'(i)$ for all $i\in\un\cN\cap \gM.$ Then $\xi$ is as 
required. 
\qed 

\subsubsection*{The ultrapower}

Let $\gN=\ult_{\cL}\,\cF$ be the ultrapower. Thus we define:\its
\bit
\item $f\approx g$ iff $\ans{i:f(i)=g(i)}\in\cL$ for $f,\,g\in 
\cF$;\its

\item $[f]=\ans{g:g\approx f}$ (the \dd\cL {\it degree\/} 
of $f$);\its

\item $[f]\hspace{2.5pt}\ain \hspace{2.5pt}[g]$ iff 
$\ans{i:f(i)\in g(i)}\in\cL$;\its

\item $\gN=\ans{[f]:f\in\cF},$ equipped with the above defined 
membership $\ain$. 
\eit

\bte
\label{1994}
$\gN$ is an elementary extension of\/ $\gM$ via the embedding 
which associates\/ $x^\ast=[\un\cN\ti\ans{x}]$ with any\/ 
$x\in\gM.$ Moreover $\gN$ is wellfounded and the ordinals in 
$\gM$ are isomorphic to the \dd\gM ordinals via the mentioned 
embedding.  
\ete
\proof See {\sc Kanovei} and {\sc van Lambalgen} [1994].
\vspace{4mm}\qed

{\it Comment\/}. 
Propositions \ref{promes} and \ref{procho} are used to 
prove the \los\/ theorem and the property of elementary 
embedding. Proposition~\ref{proreg} is used to prove 
the wellfoundedness part of the theorem. 

\section{The nature of the ultrapower}

Theorem~\ref{1994} allows to collapse $\gN$ down to a 
transitive model $\wh\gN;$ actually $\wh\gN=\ans{\wh X:X\in\gN}$ 
where 
\dm
\wh X=\ans{\wh Y:Y\in\gN\;\;\hbox{and}\;\;Y\ain X}\,.
\dm
The content of this section will be to investigate the relations 
between $\gM,$ the initial model, and $\wh\gN,$ the (transitive 
form of its) Spector ultrapower. In particular it is interesting 
how the superposition of the ``asterisk'' and ``hat'' transforms 
embeds $\gM$ into $\wh\gN$. 

\ble
\label{ee}
$x\,\longmapsto\,\wh{x^\ast}$ is an elementary embedding\/ $\gM$ 
into\/ $\wh\gN,$ equal to identity on ordinals and sets of 
ordinals (in particular on reals).
\ele
\proof Follows from what is said above.\qed\vspace{4mm}

Thus $\wh\gN$ contains all reals in $\gM.$ We now show that 
$\wh\gN$ also contains some new reals. We recall that $\ra\in\cF$ 
is a function satisfying $\ra(i)\in i$ for all $i\in\un\cN\cap\gM.$ 

Let $\Ba=\wh{[\ra]}.$ Notice that by \los\/ $[\ra]$ is a real in 
$\gN,$ therefore $\Ba$ is a real in $\wh\gN$.

\ble
\label{randa}
$\Ba$ is random over $\gM$.
\ele
\proof Let $B\sq\cN$ be a Borel set of null measure coded in $\gM;$ 
we prove that $\Ba\not\in B.$ Being of measure 0 is an absolute 
notion for Borel sets, therefore $B\cap \gM$ is a null set in $\gM$ 
as well. Corollary~\ref{mes} implies that for \dd\cL almost all  
$i,$ we have $\ra(i)\not\in B.$ By \los\/, 
$\neg\;([\ra]\hspace{2.5pt}{{\ain}}\hspace{1pt} B^\ast)$ 
in $\gN.$ Then $\Ba\not\in \wh{B^\ast}$ in $\wh\gN.$ However, by the 
absoluteness of the Borel coding, $\wh{B^\ast}=B\cap\wh\gN,$ as 
required.\qed\vspace{4mm}

Thus $\wh\gN$ contains a new real number $\Ba.$  It so 
happens that this $\Ba$ generates all reals in $\wh\gN$.

\ble
\label{lee}
The reals of\/ $\wh\gN$ are exactly the reals of\/ $\gM[\Ba]$.
\ele
\proof It follows from the known properties of random extensions 
that every real in $\gM[\Ba]$ can be obtained as $F(\Ba)$ where 
$F$ is a Borel function coded in $\gM.$ Since $\Ba$ and all reals 
in $\gM$ belong to $\wh\gN,$ we have the inclusion $\supseteq$ in 
the lemma. 

To prove the opposite inclusion let $\ba\in\wh\gN\cap\cN.$ Then 
by definition $\ba=\wh{[F]},$ where $F\in\cF.$ In turn $F=f\circ 
\ra,$ where $f\in\gM$ is a function defined on $\cN\cap\gM.$ We 
may assume that in $\gM$ $f$ maps reals into reals. Then, 
first, by Property~\ref{sol2}, $f$ is a.e. equal in $\gM$ to a 
Borel function $g=B_\ga$ where $\ga\in\cN\cap\gM$ and $B_\ga$ 
denotes, in the usual manner, the Borel subset (of $\cN^2$ in 
this case) coded by $\ga.$ Corollary~\ref{mes} shows that we have 
$F(i)=B_\ga(\ra(i))$ for \dd{\cL}almost all $i.$ In other words, 
$F(i)=B_{\ga^\ast(i)}(\ra(i))$ for \dd{\cL}almost all $i.$ 
By \los\/, this implies $[F]=B_{[\ga^\ast]}([\ra])$ in $\gN,$ 
therefore $\ba=B_\ga(\Ba)$ in $\wh{\gN}.$ By the absoluteness 
of Borel coding, we have $\ba\in\bbL[\ga,\Ba],$ therefore 
$\ba\in\gM[\Ba]$.\qed\vspace{4mm}

We finally can state and prove the principal result. 

\bte
\label{osn}
$\wh\gN\sq\gM[\Ba]$ and $\wh\gN$ coincides with\/ $\bbL^{\gM[\Ba]}
(\hbox{\rm reals}),$ the smallest subclass of\/ $\gM[\Ba]$ containing 
all ordinals and all reals of\/ $\gM[\Ba]$ and satisfying all the 
axioms of ${\bf ZF}$.
\ete
\proof Very elementary. Since $\bbV=\bbL(\hbox{reals})$ is true 
in $\gM,$ the initial Solovay model, this must be true in 
$\wh\gN$ as well. The previous lemma completes the proof.\qed

\bcor
\label{oldreals}
The set\/ $\cN\cap \gM$ of all ``old'' reals does not belong 
to $\wh\gN$.
\ecor
\proof The set in question is known to be non--measurable 
in the random extension $\gM[\Ba];$ thus it would be  
non--measurable in $\wh\gN$ as well. However $\wh\gN$ is 
an elementary extension of $\gM,$ hence it is true in 
$\wh\gN$ that all sets are measurable.\qed

\section*{References}

\ben
\item {\sc V.~Kanovei} and {\sc M.~van~Lambalgen} [1994] 
{\it Another construction of choiceless ultrapower\/}. University 
of Amsterdam, Preprint X--94--02, May 1994.

\item {\sc M.~van~Lambalgen} [1994] Independence, randomness, 
and the axiom of choice. 
{\it J. Symbolic Logic\/}, 1992, 57, 1274 -- 1304. 

\item {\sc R.~M.~Solovay} [1970] A model of set theory in which 
every set of reals is Lebesgue measurable. {\it Ann. of Math.\/}, 
1970, 92, 1 -- 56.

\item {\sc M.~Spector} [1991] Extended ultrapowers and the Vopenka -- 
Hrb\'a\v cek theorem without choice. {\it J. Symbolic Logic\/}, 1991, 
56, 592 -- 607.
\een 

\end{document}